\documentclass[a4paper,11pt]{article}

\usepackage[latin1]{inputenc}

\usepackage{amsthm,amsfonts,amsmath,amscd,amssymb}
\usepackage{a4wide}
\usepackage{array}
\usepackage{booktabs}
\usepackage{caption}
\usepackage{enumerate}
\usepackage{epstopdf}
\usepackage{float}
\usepackage{graphic x}
\usepackage{hyperref}
\usepackage{MnSymbol}
\usepackage{psfrag}
\usepackage{slashbox}
\usepackage{subfig}
\usepackage{wrapfig}

\newtheorem{proposition}{Proposition}[section]
\newtheorem{theorem}[proposition]{Theorem}
\newtheorem{lemma}[proposition]{Lemma}
\newtheorem{corollary}[proposition]{Corollary}

\newtheorem{remark}[equation]{Remark}
\newenvironment{proofof}[1]
    {\smallskip\noindent{\textbf{Proof~of~#1.}}
    \hspace{1pt}}{\hspace{-5pt}{\nobreak\quad\nobreak\hfill\nobreak
        $\square$\vspace{2pt}\par}\smallskip\goodbreak}

\newcommand{\BC}{\mathop\mathbf{BC}}                
\newcommand{\C}[1]{\mathbf{C^{#1}}}                 
\newcommand{\Lip}{\mathop\mathbf{Lip}}              
\renewcommand{\L}[1]{\mathbf{L^#1}}             
\newcommand{\W}[2]{\mathbf{W^{#1,#2}}}          
\newcommand{\naturali}{{\mathbb{N}}}                
\newcommand{\reali}{{\mathbb{R}}}                   
\newcommand{\Mp}{{\mathcal M}_{+}}              

\renewcommand{\d}[1]{\mathinner{\mathrm{d}{#1}}}    
\newcommand{\e}{\mathinner{\mathrm{e}}}         
\renewcommand{\epsilon}{\varepsilon}                
\renewcommand{\phi}{\varphi}                        
\newcommand{\modulo}[1]{{\left|#1\right|}}          
\newcommand{\norma}[1]{{\|#1\|}}                    

\numberwithin{equation}{section}
\allowdisplaybreaks[4]

\title{Splitting-Particle Methods for Structured Population Models: Convergence and Applications}

\author{J.~A.~Carrillo\thanks{Department of Mathematics,
    Imperial College London, SW7 2AZ London, United Kingdom. E-mail: {\tt
    carrillo@imperial.ac.uk}.} \and
    P.~Gwiazda\thanks{Institute of Applied Mathematics, Warsaw
    University. E-mail: \texttt{pgwiazda@mimuw.edu.pl}} \and A.~Ulikowska\thanks{Institute of Applied
    Mathematics, Warsaw University. E-mail: \texttt{aulikowska@mimuw.edu.pl}}}

\begin{document}

\maketitle

\begin{abstract}
We propose a new numerical scheme designed for a wide class of
structured population models based on the idea of operator
splitting and particle approximations. This scheme is related to
the Escalator Boxcar Train (EBT) method commonly used in biology,
which is in essence an analogue of particle methods used in
physics. Our method exploits the split-up technique, thanks to
which the transport step and the nonlocal integral terms in the
equation can be separately considered. The order of convergence of
the proposed method is obtained in the natural space of finite
nonnegative Radon measures equipped with the flat metric. This
convergence is studied even adding reconstruction and
approximation steps in the particle simulation to keep the number
of approximation particles under control. We validate our scheme
in several test cases showing the theoretical convergence error.
Finally, we use the scheme in situations in which the EBT method
does not apply showing the flexibility of this new method to cope
with the different terms in general structured population models.
\end{abstract}


\


\section{Introduction}
The main purpose of population dynamics models is to describe the
evolution of a population, which changes its size, structure, or
trait due to birth, growth, death, selection, and mutation
processes. Initially, the models are based on linear ordinary
differential equations, and as a consequence exponential growing
solutions are typically obtained. However, in many cases it is not
a realistic phenomenon, since the exponential growth can be
inhibited by environmental limitations such as lack of nutrients,
space, partners to reproduction, etc. Additionally, these models
leave out of consideration the individual's stage of development,
which strongly influences its vital functions. For example,
fertility and death rates depend heavily on the age of human
beings, the process of cell mitosis can be influenced by the age,
size or maturity level of the cell, a trait of an offspring may
depend on parents traits. Taking into account the population
structure usually leads to first order hyperbolic equations.
Finally, subsequent generations of individuals produce slight
changes in their traits due to small mutations. Selection-mutation
models typically lead to nonlocal terms due to the offspring
different trait. This paper is devoted to the numerical analysis
of such equations written in general as
\begin{align}\label{eq1}
\frac{\partial}{\partial t} \mu + \frac{\partial}{\partial
x}(b(t,\mu)\mu) + c(t,\mu) \mu & = \int_{\reali^{+}}(\eta(t,\mu)
)(y)\d \mu(y),
\\ \nonumber
\mu(0) & = \mu_o,
\end{align}
where $t \in [0,T]$ and $x \geq 0$ denote time and a structural
variable respectively, $b, c, \eta$ are vital functions depending
on $x\geq 0$, and $\mu$ is a Radon measure describing the
distribution of individuals with respect to the trait/variable
$x$. The function $b(t, \mu)$ describes the dynamics of the
transformation of the individual's state. More precisely, the
individual changes its state according to the following ODE
$$
 \dot{x} = b(t, \mu)(x).
$$
By $c(t,\mu)(x)$ we denote a rate of evolution (growth or death
rate). The integral on the right hand side accounts for an influx
of the new individuals into the system. We assume the following
form of the measure-valued function $\eta$ is of the form
\begin{equation}\label{form_of_eta}
\eta(t,\mu)(y) = \sum_{p=1}^{r}\beta_p(t,\mu)(y) \delta_{x = \bar
x_p(y)},
\end{equation}
which means that an individual at the state $y$ gives rise to
offspring being at the states $\{\bar x_p(y)\}$, $p = 1, \dots,
r$. The integral on the right-hand side has to be understood in
the B\"ochner sense, that is, by duality on test functions
$\varphi\in \mathbf{C}_0(\reali^+)$ functions as
\begin{equation}\label{form_of_eta2}
\int_{\reali^+} \int_{\reali^+}\varphi(t,x)[\d\eta(t,\mu)(y)](x)
d\mu(y) = \sum_{p=1}^{r} \int_{\reali^+} \beta_p(t,\mu)(y)
\varphi(\bar x_p(y)) \d\mu(y)\,.
\end{equation}
In case all new born individuals have the same physiological state
$x^b$, then
\begin{equation}\label{form_of_eta3}
\eta(t,\mu)(y) = \beta(t,\mu)(y) \delta_{x = x^b}\,,
\end{equation}
and the integral in \eqref{form_of_eta2} transforms into a
boundary condition. We restrict to integral operators of the form
\eqref{form_of_eta} for the sake of simplicity. In fact, the
continuous dependence of solutions of \eqref{eq1} with respect to
$\eta$ in \cite{CCC} allows for the general case to be
approximated by integral operators of the form
\eqref{form_of_eta}, and thus this restriction is done without
loss of generality, see Remark \ref{aproxeta}.

In the present paper, we develop a numerical scheme, which is
based on results obtained in \cite{Spop}, for the equation
\eqref{eq1}. It turns out that a measure setting used in the
latter paper is convenient not only from the analytical but also
from the practical numerical simulation viewpoint. Note that the
result of a measurement or an observation is usually the number of
individuals, whose state is within a specific range. For example,
demographic data provide the number of humans within certain age
cohorts. A natural way of translating such data into a
mathematical language is to make use of Dirac Deltas.

This intuitive idea was the basis for a numerical scheme called
the {Escalator Boxcar Train} (EBT) method developed in
\cite{Roos}. This method approximates in some sense a solution at
time $t$ by a sum of Dirac measures $\sum_{i\in
I}m^i(t)\delta_{x^i(t)}$. In the first step, an initial
distribution is divided into $M$ cohorts characterized by pairs
$(m^i, x^i)$, for $i=1,\dots, M$. For the $i$-th cohort, $m^i(t)$
denotes its weight at time $t$, which is a number of the
individuals within the cohort and $x^i(t)$ is its location at time
$t$, that is, an average value of the structural variable within
this cohort. The mass $m^i(t)$ changes its value due to the
process of evolution (growth or death), while $x^i(t)$ evolves
according to the characteristic lines defined by the transport
term. A boundary cohort $(m_B, x_B)$, that is, the cohort which
accounts for the influx of new individuals into the system,
evolves in a is slightly different way, since its weight changes
additionally due to the birth process. Enclosing the boundary
cohort into the system, which occurs in certain time moments, is
called the internalization process. A power of the described
method lies in its simplicity and clear biological meaning of the
output. Indeed, integrals of a population's distribution over
specified domains, which are the output, are more meaningful than
a density's value in nodal points. Originally, the EBT method was
designed for equations of the form \eqref{eq1} with the most
simplified form of the integral kernel \eqref{form_of_eta3}, and
since its invention in \cite{Roos} it has been widely used by
biologists, see e.g. \cite{EBT1,EBT2,EBT3,EBT4}.

Similar mesh-free methods called particle methods are commonly
used in problems, where one has to model a behaviour of large
groups of particles or individuals, which interact between each
other. Contrary to the EBT, particle methods were originally
designed for problems where the number of individuals was
preserved and thus the mass conservation law holds. These methods
have been successfully used for solving numerically such problems
as the Euler equation in fluid mechanics \cite{Euler2,Euler} and
Vlasov equation in plasma physics \cite{BirdLang,CR,GV}. Recently,
they are also used in problems related to crowd dynamics and
pedestrians flow \cite{PT,piccoli} or collective motion of large
groups of agents \cite{Dorsogna,kinetic,KSUB}.

As it has been stated above, in structured population models
conservation laws do not hold in general. One has to deal with new
particles, which appear due to the birth process or mutations.
Depending on the model, new individuals may appear only on the
boundary or can be distributed over the whole domain. Therefore,
one cannot exploit some natural distances for probability measures
like Wasserstein distances. The measure approach, which rigorously
deals with Dirac Deltas in models coming from biology, is
relatively new \cite{GwiazdaThomasEtAl,GwiazdaMarciniak,CCC}, and
thus a convergence of the particle based schemes for these models
was difficult to establish for a long period of time. One of the
first steps in this direction has been made for the equation
\eqref{eq1} in \cite{GwiazdaThomasEtAl,GwiazdaMarciniak}, where
existence, uniqueness, and Lipschitz dependence of solutions on
the initial data and model parameters in the space of Radon
measures were proved. By the proper choice of a metric authors
overcame the nonconservative character of the problem. Namely,
they employed a modified Wasserstein distance and the flat metric,
known also as the bounded Lipschitz distance. This framework was
the theoretical foundations for the very recent proof \cite{EBT}
of the convergence of the EBT method without any explicit error
estimates for \eqref{eq1}--\eqref{form_of_eta3}.

In this work, we shall explicitly show how the method used for
proving the well posedness of \eqref{eq1} in \cite{Spop} can be
translated into an applicable numerical scheme. We provide
estimates on the order of the convergence for the general models
\eqref{eq1}, covering in particular the case
\eqref{eq1}--\eqref{form_of_eta}. The novelty of this paper also
concerns the problem of increasing number of Dirac measures that
appears due to birth and/or mutation processes. We provide a
procedure to construct an approximation of a sum of Dirac Deltas
by a smaller amount of deltas, called the measure reconstruction
procedure, together with an error of the approximation. This paper
is organized as follows. In Section \ref{overview}, we  describe
the algorithm and the procedure of a measure reconstruction. In
Section \ref{conv_res}, we present the proof of the convergence of
the scheme together with the convergence order error analysis. In
Section \ref{simulations}, we validate our numerical scheme and
implementation by checking the convergence order in some test
cases with explicit solutions. We also use this new proposed
scheme in several examples to show the flexibility and the
accurate approximation of the evolution of the density in
structured population models even for long-time asymptotics
including cases that are not amenable for the EBT method.

\section{Particle Method}\label{overview}

\subsection{General Description}\label{GD}
The main idea of the particle method is to approximate a solution
at each time by a sum of Dirac measures. Note that even if the
initial data in \eqref{eq1} is a sum of Dirac Deltas, the integral
term possibly produces a continuous distribution at $t > 0$. This
phenomenon can be avoided due to the splitting algorithm, which
allows to separate the transport operator from the integral one
and simulate the corresponding problems successively. This is
essentially the reason why we have exploited this technique in our
scheme. To proceed with a description of the method, assume that
the approximation of the solution at time $t_k = k\Delta t$ is
provided as a sum of Dirac measures, that is,
\begin{equation}\label{initial_freeze}
\mu_{t_k}
 = \sum_{i=1}^{M_k}m^i_k\; \delta_{x^i_k},\quad M_k \in \naturali.
\end{equation}
The procedure of calculating the approximation of the solution at
time $t_{k+1}$ is divided into three main steps. In the first step
one calculates the characteristic lines for the cohorts $(m^i,
x^i)$ given by \eqref{initial_freeze}, which is equivalent to
solving the following ODE's system on a time interval $[t_k,
t_{k+1}]$:
\begin{equation}\label{b}
\frac{d}{ds} x^i(s) = b_k(x^i(s)),\quad x^i(t_k) = x_k^i,\quad i=1,\dots,M_k,
\end{equation}
where
\begin{equation}\label{b_freeze}
b_k(x) = b(t_{k}, \mu_{t_k})(x).
\end{equation}
In other words, each Dirac Delta is transported along its
characteristic to the new location $x^i_{k+1}$ without changing
its mass. The second step consists in creating new Dirac Deltas
due to the influx of new individuals and recalculating the mass of
each Dirac Delta. We have already mentioned in the introduction
above that for each $(t,\nu) \in [0,T] \times \Mp(\reali^+)$,
$\eta$ is given by
\begin{equation}\label{eta_def}
\eta(t, \nu)(y) = \sum_{p=1}^{r}\beta_p(t, \nu)(y) \; \delta_{x = \bar x_p(y)}.
\end{equation}
From this form of $\eta$, it follows that the set of possible new
states $x^l_{k+1}$ at time $t_k$ is
$$\{x^l_{k+1},\;
l=M_k+1, \dots,M_{k+1} \} :=
 \{ \bar x_p(x_{k+1}^i),\; i = 1, \dots, M_k,\; p=1,\dots, r \}.$$
Let us define
\begin{eqnarray}
\nonumber
 \mu_{k}^1 &=& \sum_{i=1}^{M_{k+1}} m^i_k\; \delta_{x^i_{k+1}},
\\[1mm]
\label{c_freeze}
c_k(x) &=& c \left(t_k,  \mu_{k}^1\right)(x),
\\[2mm]
\label{eta_freeze}
\eta_k(y) &=& \sum_{p=1}^{r}\beta_p(t_k,  \mu_{k}^1)(y) \; \delta_{x = \bar x_p(y)}
\end{eqnarray}
and for $i,j \in \{1, \dots, M_{k+1}\}$
$$
\alpha(x_{k+1}^i, x_{k+1}^j) = \left\{
   \begin{array}{l l}
     \beta_p(t_k,  \mu_{k}^1)(x_{k+1}^j), & \quad \text{if $p$ is such that } \bar x_p(x_{k+1}^j) = x_{k+1}^i,\\
     0, & \quad \text{otherwise.}
   \end{array} \right.
$$
We cannot solve an ODE system for the masses directly, since new
states will be created at any time $t_k<t<t_{k+1}$. Therefore, we
approximate it by the following explicit Euler scheme
\begin{eqnarray}\label{c_eta}
\frac{m^i_{k+1} - m^i_k}{t_{k+1} - t_k} &=& -c_k(x^i_{k+1}) m^i_k + \sum_{j=1}^{M_{k+1}} \alpha_k(x^i_{k+1},x^j_{k+1})m^j_k,
\\[1mm]
\nonumber
m^i_k &=& 0, \;\; \mathrm{for}\;\; i = M_{k}+1,\dots, M_{k+1.}
\end{eqnarray}
The resulting measure
\begin{eqnarray}\label{nu}
 \mu_{k}^2= \sum_{i=1}^{M_{k+1}}m^i_{k+1}
\delta_{x^i_{k+1}}
\end{eqnarray}
consists of $M_{k+1} \geq M_k$ Dirac Deltas. In some cases, it is
necessary to approximate the measure \eqref{nu} by a smaller
number of Dirac Deltas (see Subsection \ref{app}). If so, we
define $\mu_{t_{k+1}} =  \mathcal R (\mu_k^2) $, where $\mathcal
R(\mu_k^2)$ is the result of this approximation. Otherwise we let
$\mu_{t_{k+1}} = \mu_k^2$.

\begin{remark}\label{2approaches}
In the particular case where only one new state $x^b$ is allowed,
we can use the continuum ODE system:
\begin{eqnarray}\label{c_eta_bound} \frac{d}{ds} m^i(s) &=&
-c_k(x^i_{k+1}) m^i(s),\quad\mathrm{for}\;\; i\neq b,
\\
\nonumber
\frac{d}{ds} m^b(s) &=&  -c_k(x^b) m^b(s) + \sum_{j=1}^{M_{k+1}} \alpha_k(x^b,x^j_{k+1})m^j(s),
\end{eqnarray}
instead of the Euler approximation \eqref{c_eta}.
\end{remark}

\noindent In the method presented above, one has to deal with an
increasing number of Dirac measures, which is an important issue
to solve from the point of view of numerical simulation. In the
simplest case that all new individuals have the same size $x^b$ at
birth, then just one additional Dirac Delta is created at the
boundary at each time step. Unfortunately, in many models the
number of new particles increases so fast that after several steps
the computational cost become unacceptable. For example, in the
case of equation describing the process of cell equal mitosis, the
number of Dirac Deltas is doubled at each time step. This growth
forces us to approximate the numerical solution by a smaller
number of Dirac measures after several iterations. This procedure
is called measure reconstruction. We propose some different
methods of this reconstruction, which are discussed in the next
subsection. In order to rigorously introduce this reconstruction
procedure and to discuss the convergence of the particle method
above, we first need to introduce several distances between
measures which are relevant and useful for those purposes.

\subsection{Distances between measures}

Through this paper $\Mp(\reali^{+})$ denotes the space of
nonnegative Radon measures with bounded total variation on
$\reali^{+} = \{x \in \reali \;\colon x \geq 0\}$. We define a
metric on $\Mp(\reali^{+})$ as
\begin{equation}
  \label{distance}
  \rho_F(\mu_1, \mu_2)
  =
  \sup
  \left\{
    \int_{\reali^{+}} \phi \, \d(\mu_1-\mu_2)
    \colon
    \phi \in \C1(\reali^{+};\reali)
    \mbox{ and }
    \norma{\phi}_{\W{1}{\infty}} \leq 1
  \right\} ,
\end{equation}
where $\norma{\phi}_{\W{1}{\infty}} = \max \left\{
\norma{\phi}_{\L\infty}, \norma{\partial_x \phi}_{\L\infty}
\right\}$. \noindent $\rho_F$ is known as a \textit{flat} metric
or a \textit{bounded Lipschitz distance}. The condition
$\C1(\reali^{+}; \reali)$ in \eqref{distance} can be replaced by
$\W{1}{\infty}(\reali^{+}; \reali)$ through a standard mollifying
sequence argument applied to the test function $\phi$, as its
derivative is not involved in the value of the integral, which
implies that $\rho_F$ is the metric dual to the
$\norma{\cdot}_{({\W{1}{\infty}})^{*}}$ distance. Note that in
this paper, the space $\Mp(\reali^{+})$ is equipped with the
metric $\rho_F$ and this shall remain until said differently. The
space $(\Mp(\reali^{+}),\rho_F)$ is complete and separable.

In the following lemma we introduce $\rho$ related to $\rho_F$,
which turns out to be useful for computational purposes. Since
\cite[Theorem 6.0.2]{AmbrosioGigliSavare} gives an explicit
formula on the Wasserstein distance between two probability
measures in terms of their cumulative distribution functions, we
shall exploit this result and relate it to the flat metric. In
particular, all error estimates calculated in Section
\ref{simulations} are given in terms of $\rho$.

\begin{lemma} \label{met_eq}
Let $\mu_1, \mu_2 \in \Mp(\reali^{+})$ be such that $M_{\mu_i} =
\int_{\reali^{+}} \d \mu_i \neq 0$ and $\tilde \mu_i = \mu_i /
M_{\mu_i}$ for $i=1,2$. Define $\rho :  \Mp(\reali^{+}) \times
\Mp(\reali^{+}) \rightarrow \reali^{+}$ as the following
\begin{eqnarray}\label{rho*}
\rho(\mu_1,\mu_2) = \min\left\{M_{\mu_1}, M_{\mu_2}\right\}
W_1(\tilde \mu_1, \tilde \mu_2)+ \modulo{M_{\mu_1} - M_{\mu_2}},
\end{eqnarray}
where $W_1$ is the $1$-Wasserstein distance. Then, there exists a
constant $C_K = \frac{1}{3}\min\left\{1,
\frac{2}{\modulo{K}}\right\}$, such that
$$
C_K\rho(\mu_1, \mu_2) \leq \rho_F(\mu_1, \mu_2) \leq \rho(\mu_1,
\mu_2),
$$
where $K$ is the smallest interval such that
$\mathrm{supp}(\mu_1), \mathrm{supp}(\mu_2) \subseteq K$ and
$\modulo{K}$ is the length of the interval $K$. If $K$ is
unbounded we set $C_K = 0$.
\end{lemma}

\begin{remark}\label{rem_w}
For $\tilde \mu_1$, $\tilde \mu_2$ defined as in the lemma above,
it holds that
$$
W_1(\tilde \mu_1, \tilde \mu_2) =  \int_{0}^{1} \modulo{
F^{-1}_{\tilde \mu_1}(t) - F^{-1}_{\tilde \mu_2}(t)} \d t =
\int_{\reali^+} \modulo{F_{\tilde \mu_1}(x) - F_{\tilde \mu_2}(x)}
\d x,
$$
which follows from {\rm\cite[Section 2.2.2]{villani2}}. Since a
cumulative distribution function $F_{\mu}$ does not have to be
continuous or strictly increasing we set
$$
F^{-1}_{\mu} (s) = \sup\{x \in \reali^{+}\; : \; F_{\mu}(x) \leq
s\}, s \in [0, 1].
$$
\end{remark}

\begin{remark}
Let $\mu \in \Mp(\reali^+)$ be a probability measure and $M_1, M_2
> 0$. Then,
\begin{equation}\label{mu_diff_masses}
\rho_F(M_1 \mu, M_2 \mu) \leq \modulo{M_1 - M_2}.
\end{equation}
Indeed, let $\phi \in \C1(\reali^+; \reali)$ be such that
$\norma{\phi}_{\W{1}{\infty}} \leq 1$. Then,
\begin{eqnarray*}
\int_{\reali^+} \phi(x) \d (M_1 \mu - M_2 \mu)(x) \leq \modulo{M_1
- M_2} \int_{\reali^+} \norma{\phi}_{\L\infty} \d \mu(x) \leq
\modulo{M_1 - M_2}.
\end{eqnarray*}
Taking supremum over all admissible functions $\phi$ proves the
assertion.
\end{remark}

\begin{proofof}{Lemma \ref{met_eq}}
Let $\mu, \nu \in \Mp(\reali^{+})$ be probability measures. Assume
for the moment that $K$ is bounded, so that $\modulo{K} < +
\infty$. Note that in the definition of $W_1$
\begin{equation*}
  W_1(\mu, \nu)
  =
  \sup
  \left\{
    \int_{\reali^{+}} \phi \, \d(\mu-\nu)\;
    \colon\;
    \Lip(\phi) \leq 1
  \right\}\,,
\end{equation*}
we can assume without loss of generality that
$\norma{\phi}_{\L\infty} \leq \modulo{K}/2$. Indeed, for any
$\phi$ such that $\Lip(\phi) \leq 1$, there exists a constant $a$
and a function $\tilde \phi$ such that $\Lip(\tilde \phi) \leq 1$,
$\norma{\tilde \phi}_{\L\infty} \leq \modulo{K}/2$ and $\phi = a +
\tilde \phi$. Observe that by taking $a$ to be the middle point of
the interval $K$, and taking into account that $\Lip(\tilde \phi)
\leq 1$ and the support of the measures is included in $K$, then
$\norma{\tilde \phi}_{\L\infty} \leq \modulo{K}/2$ in $K$. Since
the values of $\phi$ can be changed arbitrarily outside $K$, then
we can assume that $\norma{\tilde \phi}_{\L\infty} \leq
\modulo{K}/2$ without loss of generality. As a consequence, we
deduce
$$
\int_{\reali^{+}}\phi(x) \d (\mu - \nu)(x) = a \int_{\reali^{+}}
\d (\mu - \nu)(x) + \int_{\reali^{+}}\tilde \phi(x) \d (\mu -
\nu)(x) = \int_{\reali^{+}}\tilde \phi(x) \d (\mu - \nu)(x),
$$
since $\int_{\reali^{+}} \d (\mu - \nu)(x)$ is equal to zero due
to the fact that $\mu$ and $\nu$ have the same mass. Therefore, we
infer that
\begin{eqnarray*}
W_1(\mu,\nu) &=& \sup
  \left\{
    \int_{\reali^{+}} \phi \, \d(\mu-\nu)\;
    \colon\;
\norma{\phi}_{\L{\infty}}\leq \modulo{K}/2,\; \Lip(\phi) \leq 1
\right\}
\\
&\leq&   \sup
  \left\{
    \int_{\reali^{+}} \phi \, \d(\mu-\nu)\;
    \colon\;
\norma{\phi}_{\W{1}{\infty}}\leq \max\{ 1, \modulo{K}/2 \}
\right\} = \max\left\{1, \frac{\modulo{K}}{2}\right\} \rho_F(\mu,
\nu).
\end{eqnarray*}
Now, let $\mu_1, \mu_2$ be as in the statement of the Lemma. Then,
\begin{eqnarray*}
\rho_F (\mu_1, \mu_2) &=& M_{\mu_1} \rho_F \left(
\frac{\mu_1}{M_{\mu_1}}, \frac{\mu_2}{M_{\mu_1}} \right) \leq
M_{\mu_1} \rho_F \left( \frac{\mu_1}{M_{\mu_1}},
\frac{\mu_2}{M_{\mu_2}} \right) + M_{\mu_1} \rho_F \left(
\frac{\mu_2}{M_{\mu_2}}, \frac{\mu_2}{M_{\mu_1}} \right)
\\
&\leq& M_{\mu_1} \rho_{F}(\tilde \mu_1, \tilde \mu_2)  +
M_{\mu_1}M_{\mu_2} \modulo{\frac{1}{M_{\mu_1}} -
\frac{1}{M_{\mu_2}}}
\\
&=& M_{\mu_1}{W_1}(\tilde \mu_1, \tilde \mu_2)  +
\modulo{M_{\mu_1} - M_{\mu_2}},
\end{eqnarray*}
where we used triangle inequality, inequality
\eqref{mu_diff_masses} and the fact that $\rho_F(\tilde \mu_1,
\tilde \mu_2) \leq W_1(\tilde \mu_1, \tilde \mu_2)$. Analogously,
we obtain
$$
\rho_F (\mu_1, \mu_2) \leq M_{\mu_2}{W_1}(\tilde \mu_1, \tilde
\mu_2)  + \modulo{M_{\mu_1} - M_{\mu_2}}
$$
and thus,
$$
\rho_F (\mu_1, \mu_2) \leq \min\left\{M_{\mu_1}, M_{\mu_2}\right\}
{W_1}(\tilde \mu_1, \tilde \mu_2)  + \modulo{M_{\mu_1} -
M_{\mu_2}} = \rho(\mu_1, \mu_2).
$$
Note that this estimate does not depend on $\modulo{K}$.

Assume that $K$ is bounded, so that the argument above applies.
Using $\phi = \pm 1$ as a test function in \eqref{distance}, we
obtain that $\modulo{M_{\mu_1} - M_{\mu_2}} \leq \rho_F(\mu_1,
\mu_2)$. Then,
\begin{eqnarray*}
\rho(\mu_1,\mu_2) &\leq& M_{\mu_1} W_1(\tilde \mu_1, \tilde \mu_2)
+ \modulo{M_{\mu_1} - M_{\mu_2}}
\\
&=& M_{\mu_1} \max\{1, \modulo{K}/2\}\rho_F(\tilde \mu_1, \tilde
\mu_2) + \modulo{M_{\mu_1} - M_{\mu_2}}
\\
&\leq&
 \max\{1, \modulo{K}/2\} \rho_F\left( \mu_1, \frac{M_{\mu_1}}{M_{\mu_2}} \mu_2\right) + \rho_F(\mu_1, \mu_2)
\\
&\leq&
 \max\{1, \modulo{K}/2\} \left(\rho_F\left( \mu_1, \mu_2 \right) +
\rho_F\left( \mu_2, \frac{M_{\mu_1}}{M_{\mu_2}}
\mu_2\right)\right) +
 \rho_F(\mu_1, \mu_2)
\\
&\leq& 2\max\{1, \modulo{K}/2\} \rho_F(\mu_1, \mu_2) + \max\{1,
\modulo{K}/2\} M_{\mu_2}  \modulo{1 - \frac{M_{\mu_1}}{M_{\mu_2}}}
\\
&\leq& 3\max\{1, \modulo{K}/2\} \rho_F(\mu_1, \mu_2),
\end{eqnarray*}
which implies that
$$
\frac{1}{3}\min\left\{1, \frac{2}{\modulo{K}}\right\}
\rho(\mu_1,\mu_2) \leq \rho_F(\mu_1,\mu_2).
$$
In case $\modulo{K} = +\infty$ we set $C_K = 0$ obtaining a
trivial inequality $0 \leq \rho_F(\mu_1, \mu_2)$.
\end{proofof}
\begin{remark}
The dependence of the constant $C_K$ on a length of the interval
$K$ express a small sensitivity of the flat metric in the case
where a distance between supports of measures is large. In
particular, the flat distance for two Dirac measures $\delta_{x =
a}$ and $\delta_{x=b}$ is equal to
$$
\rho_F(\delta_{x=a}, \delta_{x=b}) = \min\{2, \modulo{a-b}\}.
$$
\end{remark}

Now, we can precisely discuss the measure reconstruction by
approximation with a fixed number of particles of continuum or
larger number of particles distributions.

\subsection{Measure Reconstruction}\label{app}
Due to Lemma \ref{met_eq}, we restrict our analysis to
probability measures. Let $\mu = \sum_{i=1}^{M} m_i \delta_{x_i}$
be a probability measure with a compact support $K=[k_1,k_2]$. The
aim of the reconstruction is to find a smaller number of Dirac
Deltas $\bar M<M$ such that
\begin{eqnarray*}
\mathcal R_o(\mu) := \mathrm{argmin}\;\; W_1 \left(\mu,
\sum_{j=1}^{\bar M} n_j \delta_{y_j} \right),\quad
\mathrm{where}\;\; \sum_{j=1}^{\bar M} n_j = 1
\;\;\mathrm{and}\;\; n_j \geq 0, x_j \in \reali^+.
\end{eqnarray*}
This minimisation procedure is essentially a linear programming
problem which, under some particular assumptions on cycles, can be
solved by the simplex algorithm providing the global minimum. This
choice is the optimal for the reconstruction procedure. However,
its complexity is at least cubic. From that reason, we exploit
less costly (linear cost in the size of the problem) methods of
reconstruction, which provide the error of the order $\mathcal
O(1/\bar M)$. Note that the cubic cost is unacceptable in our
case, since the total cost of the method is quadratic if the
number of particles grows linearly with the time step.

\

\noindent \textbf{A) Fixed-location reconstruction:} The idea of
the fixed-location reconstruction is to divide the   support of
the measure $\mu$ into $\bar M$ equal intervals and put a Dirac
Delta with a proper mass in the middle of each interval. The mass
of this Dirac Delta is equal to the mass of $\mu$ contained in
this particular interval. Let $\Delta x =  \modulo{K}/\bar M$ and
define
$$
\tilde{x}_j = k_1 + \left(j - \frac{1}{2} \right)\Delta x, \quad
\tilde m_j  = \left\{
\begin{array}{l l}
\mu\left(\left[\tilde{x}_j - \Delta x/2,\tilde{x}_j +\Delta x/2\right)\right),\;
&
\text{for}\; j=1,\dots, \bar M-1,
\\[1mm]
 \mu\left(\left[\tilde{x}_{\bar M} - \Delta x/2,\tilde x_{\bar M} +\Delta x/2\right]\right),\;
& \text{for}\; j=\bar M,
\end{array}
\right.
$$
and
$$
\mathcal R_l(\mu) := \sum_{j=1}^{\bar M} \tilde m_j
\delta_{\tilde{x}_j}\,.
$$
To estimate the error between $\mu$ and $\mathcal R_l(\mu)$
consider a transportation plan $\gamma$ between both measures.
Then, according to \cite[Introduction]{villani2}, we have
\begin{eqnarray}\label{e1}
{W_1} \left(\mu, \mathcal R_l(\mu) \right) \leq
\int_{\reali_{+}^2} \modulo{\; x - y \;} \d \gamma(x,y) \leq
\int_{\reali_{+}^2} \frac{\Delta x}{2}   \d \gamma(x,y) \leq
\frac{\Delta x}{2} = \frac{\modulo{K} }{ 2\bar M}.
\end{eqnarray}
The second inequality follows from the fact that each particle was
shifted by a distance not greater than a half of the interval of a
length $\Delta x$, while the third one is a consequence of the
fact that $\gamma$ is a probability measure on $\reali_{+}^{2}$.

\

\textbf{B) Fixed-Equal mass reconstruction:} The aim of the
fixed-equal mass reconstruction is to distribute Dirac Deltas of
equal masses over the support of a given measure in a proper way.
In our particular case, we want to reduce the number of Dirac
Deltas from $M$ to $\bar M$, and thus we need to explain an
algorithm allowing for splitting of the Dirac Deltas into two. The
definition of the reconstruction operator $\mathcal R_m(\mu)$ is
as follows: we set
$$
\tilde{m}_j = \frac1{\bar M},\;\;\mathrm{for}\;\; j=1,\dots,\bar
M.
$$
The scheme for determining $\tilde{x}_j$ is the following. We
first look for an index $n_1$, such that
$$
\sum_{i=1}^{n_1 - 1}m_i < \frac1{\bar M} \leq \sum_{i=1}^{n_1}m_i.
$$
We set
$$
\tilde{x}_1 = \sum_{i=1}^{n_1 - 1} m_i x_i + m'_{n_1} x_{n_1}, \;\;\mathrm{where}\;\;
m'_{n_1} = \frac1{\bar M} - \sum_{i=1}^{n_1 - 1} m_i x_i.
$$
Namely, the mass located in $x_{n_1}$ is split into two parts --
the amount of mass equal to $m'_{n_1}$ is shifted to $\tilde{x}_1$
and the rest, that is, $m_{n_1} - m'_{n_1}$ stays in $x_{n_1}$.
For simplicity, we redefine $m_{n_1} : = m_{n_1} - m'_{n_1}$ and
repeat the procedure described above until the last point
$\tilde{x}_{\bar M}$ is found to get the final form of the
reconstruction
$$
\mathcal R_m(\mu) := \sum_{j=1}^{\bar M} \tilde m_j
\delta_{\tilde{x}_j}\,.
$$

Note that in each step of the procedure one changes the locations
of the Dirac Deltas, of which joint mass is not greater than $m$.
Using an analogous argument as in the previous case, we conclude
that in the $j$-th step we commit an error not greater than
$\modulo{x_{n_j} - x_{n_{j-1}}}m$, where $x_{n_o} = k_1$. Since
$k_1 = x_{n_{o}} \leq x_{n_{1}}\dots \leq x_{n_{\bar M}} \leq
k_2$, the total error can be bounded by
\begin{eqnarray}\label{e2}
W_1(\mu, \mathcal R_m(\mu) ) \leq \frac{\modulo{K}}{\bar M}.
\end{eqnarray}

\

The findings above can be summarized in the following

\begin{corollary}\label{correcon}
The error of the fixed-location $\mathcal R_l(\mu)$ and
fixed-equal mass $\mathcal R_m(\mu)$ reconstructions is of the
order of $\mathcal{O}(1/\bar M)$ where $\bar M$ is the number of
Dirac Deltas approximating the original measure $\mu$.
\end{corollary}

These reconstructions can be used at $t=0$, if the initial data in
\eqref{eq1} is not a sum of Dirac Deltas or at $t>0$ in order to
deal with the problem of increasing number of Dirac Deltas, which
are produced due to birth and/or mutation processes. We introduce
the following notation:
\begin{itemize}
\item $E_I(\bar M_o)$ is the upper bound for the error of the
initial data reconstruction defined in terms of $W_1$ distance.
More specifically, for a measure $\mu$ such that $M_{\mu} :=
\int_{\reali^+} \d \mu(x)> 0$, it holds that
$$
W_1\left( \frac{\mu}{M_{\mu}}, \frac{\mathcal R(\mu)}{M_{\mu}}\right) \leq E_I(\bar M_o).
$$
Here, the reconstruction operator $\mathcal R(\mu)$ refers to
either $\mathcal R_l(\mu)$ or $\mathcal R_m(\mu)$.

\item $E_R(\bar M)$  is the upper bound for the error of the
measure reconstruction at time $t > 0$ defined in terms of $W_1$
distance as above.
\end{itemize}

We are now ready to state and prove the main convergence result.


\section{Convergence Results}\label{conv_res}

\subsection{Assumptions and theoretical results on splitting}
For the sake of the reader, we recall the theoretical results on
splitting for the equation \eqref{eq1} obtained in \cite{Spop}.
The assumptions on the parameter functions $b,c$ and $\beta_p$,
$p=1, \dots, r$, are the following
\begin{eqnarray}
  \label{eqAssumptions}
    b,c, \beta_p \; : \;
    [0, T] \times \mathcal{M}^+(\reali^+)
    &\to&
    \W{1}{\infty}(\reali^+;\reali),
\\
\label{eqAssumptions_xp}
\bar x_p \; : \; \reali^+ &\to&
\reali^+,
  \end{eqnarray}
where $ b(t,\mu)(0)\geq 0$ for $(t,\mu)\in [0,T] \times
\mathcal{M}^+(\reali^+)$ and $p=1, \dots, r$. We require the
following regularity
  \begin{eqnarray}\label{assumptions:nonlinear}
  b,c, \beta_p
  & \in &
  {\BC}^{\mathbf{\alpha,1}}
  \left(
    [0,T] \times {\mathcal M}^{+}(\reali^+); \;
    \W{1}{\infty}(\reali^+;\reali)
  \right),
\\
\label{assumptions:nonlinear_xp}
\bar x_p &\in &
\Lip(\reali^+; \reali^+).
\end{eqnarray}
Here, ${\BC}^{\mathbf{\alpha,1}}([0,T]\times{\mathcal
M}^+(\reali^+); \W{1}{\infty}(\reali^+;\reali))$ is the space of
$\W{1}{\infty}(\reali^+;\reali)$ valued functions which are
bounded in the $\norma{\cdot}_{\W{1}{\infty}}$ norm, H\"older
continuous with exponent $0<\alpha\leq 1$ with respect to time and
Lipschitz continuous in $\rho_F$ with respect to the measure
variable. This space is equipped with the
$\norma{\cdot}_{{\BC}^{\mathbf{\alpha,1}}}$ norm defined by
\begin{equation}
\label{norma_bc}
  \norma{f}_{\mathbf{\BC^{\alpha,1}}}
  =
  \sup_{t\in[0,T], \mu\in{\mathcal{M}^+(\reali^+)}}
  \left(
\norma{f(t, \mu)}_{\W{1}{\infty}} +
    \Lip\left(f(t,\cdot)\right) +
    \mathrm{H}_\alpha\left(f(\cdot,\mu)\right)
  \right),
\end{equation}
where $\Lip(f)$ is the Lipschitz constant of a function $f$ and
\begin{displaymath}
  \mathrm{H}_\alpha(f(\cdot,\mu)) :=
  \sup_{s_1,s_2\in[0,T]}
    \frac{\norma{f(s_1,\mu) - f(s_2,\mu)}_{\W{1}{\infty}}}{\modulo{s_1 -
    s_2}^{\alpha}} .
\end{displaymath}
For any $f \in {\BC}^{\mathbf{\alpha,1}}([0,T]\times{\mathcal
M}^+(\reali^+);\W{1}{\infty}(\reali^+;\reali))$ and any $\mu:[0,T]
\to \mathcal M^+(\reali^+)$, we define
\begin{equation*}
\norma{f}_{\BC} = \sup_{t\in[0,T]}\norma{f(t, \mu(t))}_{\L\infty}.
\end{equation*}
Regularity of  $\beta_p$ and $x_p$ imposed in
\eqref{eqAssumptions}--\eqref{assumptions:nonlinear_xp} guarantees
that $\eta$ defined by \eqref{form_of_eta} fulfills the
assumptions of \cite[Theorem 2.11]{Spop} and thus, \eqref{eq1} is
well posed. We recall this result next.

\begin{theorem}\label{thm:Main}
Let \eqref{eqAssumptions}--\eqref{assumptions:nonlinear_xp} hold.
Then, there exists a unique solution
$$
  \mu \in (\BC\cap
  \Lip)\left([0,T] ;{\mathcal M}^+(\reali^+)\right)$$ to~\eqref{eq1}.
  Moreover, the following properties are satisfied:
  \begin{enumerate}
  \item For all $0 \leq t_1\leq t_2 \leq T$ there exist constants
    $K_1$ and $K_2$, such that
    \begin{equation*}
      \rho_F\left(\mu(t_1),\mu(t_2)\right)
      \leq
      K_1 \e^{K_2({t_2-t_1})} \mu_o(\reali^+)({t_2 - t_1}).
    \end{equation*}
  \item Let $\mu_1(0), \mu_2(0) \in {\mathcal M}^{+}(\reali^+)$ and
    $b_i$, $c_i$, $\beta_i = (\beta^i_1, \dots, \beta^i_r)$ satisfy
    assumptions \eqref{eqAssumptions} - \eqref{assumptions:nonlinear_xp} for $i = 1, 2$, $p=1,\dots, r$. Let $\mu_i$
    solve~\eqref{eq1} with initial datum $\mu_i(0)$ and
    coefficients $(b_i,c_i,\beta_i)$.  Then, there exist
    constants $C_1$, $C_2$ and $C_3$ such that for all $t\in [0,T]$
    \begin{eqnarray*}
      \rho_F\left(\mu_1(t),\mu_2(t)\right)
      \leq \e^{C_1t} \rho_F\left(\mu_1(0),\mu_2(0)\right)
      +
      C_2 \e^{C_3t}t\;
      \norma{(b_1,c_1,\beta_1) - (b_2,c_2, \beta_2)}_{\BC}.
    \end{eqnarray*}
    where
    $$
\norma{(b,c,\beta)}_{\BC} = \norma{b}_{\BC} + \norma{c}_{\BC} +
\sum_{p=1}^{r}\norma{\beta_p}_{\BC} \,.
    $$
  \end{enumerate}
\end{theorem}

\subsection{Error estimates in $\rho_F$}
The aim of this subsection is to obtain an estimate on the error
between the numerical solution $\mu_t$ and the exact solution
$\mu(t)$. Let $[0,T]$ be a time interval, $N$ be a number of time
steps, $\Delta t = T / N$ be the time step. We define the time
mesh $\{t_k\}_{k=0}^N$, where $t_k = k \Delta t$. Let $\bar M_k$,
$k = 0,1,\dots N$, be parameters of the measure reconstruction. In
particular, $\bar M_o$ is the number of Dirac Deltas approximating
the initial condition and $\bar M_k$ stands for the number of
Dirac measures approximating the numerical solution at $t > 0$
after a reconstruction, if performed. We assume that
reconstructions are done every $n$ steps, which means that there
are $\mathcal K = N/n$ reconstructions, each at time $t_{jn}$,
where $j = 1,\dots,\mathcal K$. Let $\bar M$ be the number of
Dirac Deltas after the reconstruction that will not depend on
time.
\begin{theorem}\label{rec_error}
Let $\mu$ be a solution to \eqref{eq1} with initial data $\mu_o$.
Assume that $\mu_{t_m}$ is defined by the numerical scheme
described in Subsection {\rm \ref{GD}} and $m = jn$ for some $j
\in \{1,\dots, \mathcal K\}$, i.e., that $t_m$ is the time after
$j$ reconstructions. Then, there exists $C$ depending only on the
parameter functions, the initial data, and $T$ such that
\begin{equation}\label{ppp}
\rho_F\left( \mu_{t_m},\mu(t_m) \right) \leq C\left( \Delta t +
(\Delta t)^{\alpha} + E_I(\bar M_o) + E_R(\bar M) j\right).
\end{equation}
\end{theorem}

\begin{remark}\label{rmerror}
The error estimate \eqref{ppp} accounts for different error
sources. More specifically, the error of the order
$\mathcal{O}(\Delta t)$ is a consequence of the splitting
algorithm. The term of order $\mathcal{O}((\Delta t)^{\alpha})$
follows from the fact that we solve \eqref{b}--\eqref{c_eta} with
parameter functions independent of time, while $b,c$ and $\eta$
are in fact of $\C{\alpha}$ regularity with respect to time.
Finally, $E_I$ and $E_R$ are the errors coming from the measure
reconstruction procedure that are of the order $1/\bar{M_o}$ and
$1/\bar{M}$ respectively as proven in subsection {\rm 2.3}.
Thinking about $1/\bar{M}$, with $\bar{M}=\bar{M}_o$, as the
spatial discretization $\Delta x$ and for $\alpha=1$, we obtain
that the method is of order one both in space and in time.
\end{remark}
\begin{proofof}{ Theorem \ref{rec_error}}
The proof is divided into several steps. For simplicity, in all
estimates below, we will use a generic constant $C$, without
specifying its exact form that may change from line to line.
\smallskip
\\
\textbf{Step 1: The auxiliary scheme.} \quad Let us define the
auxiliary semi-continuous scheme, which consists in solving
subsequently the following problems:
\begin{equation}
\left\{
\begin{array}{rcl}
\label{exact1}
\partial_t \mu + \partial_x ( \bar b_k(x)\mu )&=& 0,\quad
\mathrm{on}\;\;\; [t_k, t_{k+1}] \times \reali^+,
\\
\mu(t_{k}) &=& \bar \mu_k
\end{array}
\right.
\end{equation}
and
\begin{equation}
\label{exact2}
\left\{
\begin{array}{rcl}
\partial_t \mu &=& -\bar{\bar c}_k(x)\mu + \int_{\reali^+} \bar{ \bar{\eta}}_k(y) \d \mu(y),\quad
\mathrm{on}\;\;\; [t_k, t_{k+1}] \times \reali^+,
\\
\displaystyle
\mu(t_{k}) &=& \bar \mu_k^1,
\end{array}
\right.
\end{equation}
where $\bar \mu_k \in \mathcal M^+(\reali^+)$, $\bar \mu_k^1$ is
the solution to \eqref{exact1} at time $t_{k+1}$ and $\bar b_k$,
${\bar{\bar c}}_k$, and $\bar{\bar \eta}_k$ are defined as
\begin{eqnarray}
\label{bar_b_freeze}
\bar b_k(x) &=& b \left(t_k,  \bar \mu_{k}\right)(x),
\\[2mm]
\nonumber
{\bar{\bar c}}_k(x) &=& c \left(t_k,  \bar \mu_{k}^1\right),
\quad
\bar{\bar {\eta}}_k(y) = \sum_{p=1}^{r}\beta_p(t_k,  \bar \mu_{k}^1)(y) \; \delta_{x = \bar x_p(y)}.
\end{eqnarray}
A solution to the second equation at time $t_{k+1}$ is denoted by
$\bar \mu^2_k$. The output of one time step of our scheme is
defined as $\bar \mu_{k+1} =\mathcal R (\bar \mu_k^2)$.
\smallskip
\\
\textbf{Step 2: Error of the reconstruction.} \quad Since $\bar
\mu_{k+1}$ arises from $\bar \mu_k^2$ through the reconstruction,
masses of both measures are equal. Therefore, application of Lemma
\ref{met_eq} yields
\begin{eqnarray}\label{z1}
\rho_F(\bar \mu_{k+1},\bar \mu_k^2) \leq \rho(\bar \mu_{k+1}, \bar \mu_k^2)
=
M_{\bar \mu_k^2}
  W_1 \left( \frac{\bar \mu_{k+1}}{M_{\bar \mu_k^2}}, \frac{\bar \mu_k^2}{M_{\bar \mu_k^2}}\right) \leq M_{\bar \mu_k^2}  E_R(\bar M),
\end{eqnarray}
where $M_{\bar \mu_k^2} =\bar \mu_{k+1}(\reali^{+}) = \bar
\mu_k^2(\reali^{+})$ and $E_R(\bar M)$ is the error of the
reconstruction introduced in Subsection \ref{app}. As stated in
Corollary \ref{correcon}, $E_R(\bar M)$ is of order $1/\bar M$ for
both reconstructions. Note that $M_{\bar \mu_k^2} $ can be bounded
independently on $k$. Indeed, on each time interval $[t_k,
t_{k+1}]$ mass grows at most exponentially, which follows from
\cite[Theorem 2.10, (i)]{Spop}, and reconstructions, if performed,
do not change the mass. Thus, there exists a constant $C = C(T, b,
c, \eta, \mu_o)$ such that $M_{\bar \mu_k^2}  \leq C$.

\textbf{Step 3: Error of splitting.} \quad Let $\nu(t)$ be a
solution to \eqref{eq1} on a time interval $[t_{k}, t_{k+1}]$ with
initial datum $\bar \mu_{k}$ and parameter functions $\bar b_k$,
$\bar c_k$, $\bar \eta_k$, where $\bar b_k$ is defined by
\eqref{bar_b_freeze},
\begin{eqnarray}
\label{bar_c_freeze}
{{\bar c}}_k(x) &=& c \left(t_k,  \bar \mu_{k}\right),
\\
\label{bar_eta_freeze}
\bar{ {\eta}}_k(y) &=& \sum_{p=1}^{r}\bar{\beta}_{p,k}(y) \; \delta_{x = \bar x_p(y)},\quad
\mathrm{where}\;\;\; \bar{\beta}_{p,k}(y) = \beta_p(t_k,  \bar \mu_{k})(y).
\end{eqnarray}
According to \cite[Proposition 2.7]{ColomboGuerra2009} and
\cite[Proposition 2.7]{Spop}, the distance between $\bar
\mu_{k}^2$ and $\nu(t_{k+1})$, that is, the error coming from the
application of the splitting algorithm can be estimated as
\begin{equation}\label{est_split}
\rho_F(\bar \mu_{k}^2,\nu(t_{k+1})) \leq C (\Delta t)^2.
\end{equation}
To estimate a distance between $\nu(t_{k+1})$ and $\mu(t_{k+1})$
consider $\zeta(t)$, which is a solution to \eqref{eq1} on a time
interval $[t_k, t_{k+1}]$ with initial data $\mu(t_k)$ and
coefficients $\bar b_k$, $\bar c_k$, $\bar \eta_k$. By triangle
inequality
$$
\rho_F(\nu(t_{k+1}), \mu(t_{k+1})) \leq \rho_F(\nu(t_{k+1}), \zeta(t_{k+1})) + \rho_F(\zeta(t_{k+1}), \mu(t_{k+1})).
$$
The first term of the inequality above is a distance between
solutions to \eqref{eq1} with different initial data, that is,
$\bar \mu_{k}$ and $\mu(t_k)$ respectively. The second term is
equal to a distance between solutions to \eqref{eq1} with
coefficients $(\bar b_k, \bar c_k,\bar \eta_k)$ defined by
\eqref{bar_b_freeze}, \eqref{bar_c_freeze}, and
\eqref{bar_eta_freeze}, and $(b(t,\mu(t)), c(t,\mu(t)),
\eta(t,\mu(t)))$. By the continuity of solutions to \eqref{eq1}
with respect to the initial datum and coefficients in Theorem
\ref{thm:Main}, we obtain
\begin{eqnarray}\label{main1}
\rho_F(\nu(t_{k+1}), \zeta(t_{k+1}))
\leq
\e^{ C \Delta t} \rho_F(\bar \mu_{k}, \mu(t_{k})),
\end{eqnarray}
and
\begin{eqnarray}\label{main2}
\rho_F(\zeta(t_{k+1}), \mu(t_{k+1})) \leq C \Delta t  \e^{ C\Delta
t} \left( \norma{\bar b_k - b}_{\overline{\BC}} + \norma{\bar c_k
- c}_{\overline{\BC}} + \sum_{p=1}^{r} \norma{\bar \beta_{p,k} -
{\beta_p}}_{\overline{\BC}} \right),
\end{eqnarray}
where
\begin{eqnarray}
\nonumber \norma{\bar b_k - b}_{\overline{\BC}} &=& \sup_{t \in
[t_k, t_{k+1}]}\norma{ \bar b_k - b(t,\mu(t))}_{\L\infty},
\\
\label{111} \norma{\bar c_k - c}_{\overline{\BC}} &=& \sup_{t \in
[t_k, t_{k+1}] }\norma{\bar c_k - c(t,\mu(t))}_{\L\infty},
\\
\label{222} \norma{\bar \beta_{p,k} - \beta_p}_{\overline{\BC}}
&=& \sup_{t \in [t_k, t_{k+1}]} \norma{\bar \beta_{p,k} -
\beta_p(t, \mu(t))}_{\L\infty}.
\end{eqnarray}
Due to the assumptions
\eqref{eqAssumptions}--\eqref{assumptions:nonlinear_xp} and the
definition of $\bar b_k$, $\bar c_k$, $\bar \eta_k$ we obtain
\begin{eqnarray}
\nonumber
\norma{\bar b_k - b(t,\mu(t))}_{\L\infty}
&\leq&
\norma{ b(t_{k},\bar \mu_{k}) - b(t_{k},\mu(t))}_{\L\infty} +
\norma{b(t_{k},\mu(t)) - b(t,\mu(t))}_{\L\infty}
\\
\nonumber
&\leq&
 \Lip(b(t_k, \cdot))\; \rho_F(\bar \mu_{k}, \mu(t)) +
 \norma{b}_{\BC^{\mathbf {\alpha,1}}} \modulo{t - t_k}^{\alpha}
\\
\label{333} &\leq& \norma{b}_{\BC^{\mathbf {\alpha,1}}}\left[
\rho_F(\bar \mu_{k}, \mu(t)) + (\Delta t)^{\alpha}\right]\,.
\end{eqnarray}
Using Lipschitz continuity of the solution $\mu(t)$, see
\cite[Theorem 2.11]{Spop}, we obtain
$$
\rho_F(\bar \mu_{k}, \mu(t))\leq
 \rho_F(\bar \mu_{k}, \mu(t_k)) + \rho_F(\mu(t_k), \mu(t))
\leq
 \rho_F(\bar \mu_{k}, \mu(t_k))  + C \Delta t \e^{ C \Delta t}.
$$
Substituting the latter expression into \eqref{333} yields
$$
\norma{b_k - b(t,\mu(t))}_{\L\infty} \leq
\norma{b}_{\BC^{\alpha,1}} \left(
 \rho_F(\bar \mu_{k}, \mu(t_k))  + C \Delta t \e^{ C \Delta t}
\right)
+ \norma{b}_{\BC^{\alpha,1}}(\Delta t)^{\alpha}.
$$
Bounds for \eqref{111} and \eqref{222} can be proved analogously.
From the assumptions it holds that
$$
\norma{(b,c,\beta)}_{\BC^{\mathbf{\alpha,1}}} = \norma{b}_{\BC^{\mathbf{\alpha,1}}} + \norma{c}_{\BC^{\mathbf{\alpha,1}}} + \sum_{p=1}^{r}\norma{\beta_p}_{\BC^{\mathbf{\alpha,1}}} < +\infty,
$$
and as a consequence, we obtain
$$
\norma{\bar b_k - b}_{\overline{\BC}} + \norma{\bar c_k -
c}_{\overline{\BC}} +\sum_{p=1}^{r} \norma{\bar \beta_{p,k} -
{\beta_p}}_{\overline{\BC}} \leq
\norma{(b,c,\beta)}_{\BC^{\mathbf{\alpha,1}}} \left[ \rho_F(\bar
\mu_{k}, \mu(t_k)) + C\e^{ C  T} \Delta t + (\Delta
t)^{\alpha}\right]\,.
$$
Using this inequality in \eqref{main2} yields
\begin{eqnarray*}
\rho_F(\zeta(t_{k+1}), \mu(t_{k+1})) &\leq& C \Delta t  \e^{
C\Delta t} \left[
 \rho_F(\bar \mu_{k}, \mu(t_k))
+
 \Delta t
+
 (\Delta t)^{\alpha} \right]
\\
&\leq&
C \Delta t  \e^{ C\Delta t}
 \rho_F(\bar \mu_{k}, \mu(t_k))
+ C   \e^{ C T} (\Delta t)^2 + C \e^{ C T} (\Delta
t)^{1+\alpha}\,.
\end{eqnarray*}
Combining the inequality above with \eqref{main1} and redefining
$C$ leads to
\begin{eqnarray}\label{z2}
\rho_F(\nu(t_{k+1}), \mu(t_{k+1})) &\leq& \e^{ C \Delta t}(1 +
C\Delta t) \rho_F(\bar \mu_{k}, \mu(t_{k})) + C (\Delta t)^2 + C
(\Delta t)^{1+\alpha}
\\[2mm]
\nonumber &\leq& \e^{ 2C \Delta t} \rho_F(\bar \mu_{k},
\mu(t_{k})) + C (\Delta t)^2 + C (\Delta t)^{1+\alpha}.
\end{eqnarray}
Finally, putting together \eqref{z2} and \eqref{est_split}, we
conclude that
\begin{equation}\label{znew}
\rho_F(\bar \mu_k^2, \mu(t_{k+1})) \leq \e^{ 2C \Delta t}
\rho_F(\bar \mu_{k}, \mu(t_{k})) + C (\Delta t)^2 + C (\Delta
t)^{1+\alpha}.
\end{equation}

\textbf{Step 4: Adding the errors.} \quad Now, let $w = j n$, $v =
(j-1)n$, $j \in \{1,\dots,\mathcal K\}$, that is, $t_w$ and
$t_{v}$ are the time points in which the measure reconstruction
occurs. Since for $t_i$ such that $t_v < t_i < t_w$ it holds that
$\bar \mu_{i} = \mathcal R (\bar \mu_{i-1}^2) = \bar \mu_{i-1}^2$,
i.e., the measure reconstruction is not performed, the application
of the discrete Gronwall's inequality to \eqref{znew} yields
\begin{eqnarray*}
 \rho_F(\bar \mu_w^2, \mu(t_w))
&\leq&
 \e^{nC\Delta t}  \rho_F(\bar \mu_v, \mu(t_v))
+ C \frac{\e^{n C \Delta t}-1}{\e^{C\Delta t} - 1}  \left( (\Delta
t)^2 + (\Delta t)^{1+\alpha}\right).
\end{eqnarray*}
There exists $C^*$ depending only on $T$ such that $\e^{n C \Delta
t}-1 < n C^* \Delta t$, for each $n \Delta t \in [0,T]$.
Therefore, we deduce
$$
\frac{\e^{nC \Delta t}-1}{\e^{C\Delta t} - 1} \leq
\frac{nC^*\Delta t}{C \Delta t} = \frac{C^*}{C}n
$$
and thus,
\begin{equation*}
 \rho_F(\bar \mu_w^2, \mu(t_w))
\leq
 \e^{nC\Delta t}  \rho_F(\bar \mu_v, \mu(t_v))
+ n C \left( (\Delta t)^2 + (\Delta t)^{1+\alpha}\right),
\end{equation*}
for some constant $C$. Combining this inequality with \eqref{z1}
in Step 2 of the proof yields
\begin{eqnarray*}
 \rho_F(\bar \mu_w, \mu(t_{w})) &\leq&
 \e^{n C \Delta t}  \rho_F(\bar \mu_v, \mu(t_v))
+ nC ((\Delta t)^2 + (\Delta t)^{1+\alpha}) + C E_R(\bar M).
\end{eqnarray*}
\textbf{Step 5: Final estimate for the auxiliary scheme.} An
analogous argument using the discrete Gronwall's inequality again
results in the following estimate
\begin{eqnarray}
\nonumber
 \rho_F(\bar \mu_w, \mu(t_{w}))
&\leq&\!\!
 \e^{jn C\Delta t}  \rho_F(\mathcal R(\mu_o), \mu_o)
+ C \frac{ \e^{jn C\Delta t} - 1}{ \e^{n C\Delta t} - 1} \left[n(
(\Delta t)^2 + (\Delta t)^{1+\alpha}) + E_R(\bar M) \right]
\\
\nonumber &\leq& C \e^{C t_w} E_I(\bar M_o) + Cj  \left[n (
(\Delta t)^2 + (\Delta t)^{1+\alpha}) + E_R(\bar M) \right]
\\
&\leq& C\e^{C t_w} E_I(\bar M_o) + C(jn \Delta t) \left(\Delta t +
(\Delta t)^{\alpha}\right) + Cj E_R(\bar M) \label{final}
\end{eqnarray}
and since $jn \Delta t = t_w \leq T$ the assertion is proved.
\smallskip
\\
\noindent \textbf{Step 6: Full error estimate.} The full error
estimate \eqref{ppp} takes into account the error following from
the numerical approximation of the auxiliary scheme
\eqref{exact1}--\eqref{exact2}. This additional source of error is
nothing else than the error of the Euler method for ODE's.
According to \cite[(515.62)]{butcher}, the error committed is of
order $\Delta t$ when solving \eqref{exact1}--\eqref{exact2} using
its Euler approximation \eqref{b}--\eqref{c_eta}. Therefore, the
final estimate \eqref{final} holds.
\end{proofof}

\begin{remark}\label{aproxeta}
In this work, we have assumed that $\eta$ is given as a sum of
Dirac Deltas \eqref{form_of_eta}. If $\eta(t,\mu)(y)$ is not in
such a form, one has to use a proper approximation by atomic
measures in order to apply our scheme. One of the possibilities
for this approximation is through the measure reconstruction
described in Subsection {\rm\ref{app}}. Assume that there exists a
bounded interval $K$ such that for all $(t,\mu) \in [0,T] \times
\mathcal M^+(\reali^+)$, we have
\begin{equation}\label{ttt}
\mathrm{supp}(\eta(t,\mu)(y)) \subseteq K.
\end{equation}
Fix $r \in \naturali$ and let $\{K_p\}_{p=1}^r$ be a family of intervals such that
$$
\bigcup_{p=1}^r K_p = K,
\quad
K_i \cap K_j = \emptyset,\; \mathrm{for}\; i \neq j
\quad\mathrm{and} \quad
\modulo{K_p} = \frac{\modulo{K}}{r},\;\mathrm{where}\; p=1,\dots, r.
$$
Namely, we divide $K$ into $r$ disjoint intervals of equal length.
Denote the center of each interval by $\bar x_p(y)$ and define
\begin{equation}\label{mmm2}
\beta_p(t,\mu)(y) = \int_{K_p} \d (\eta(t,\mu)(y)) (x).
\end{equation}
The approximation of $\eta(t,\mu)(y)$ is thus given by
\begin{equation}\label{mmm}
\sum_{p=1}^{r} \beta_p(t,\mu)(y) \delta_{x = \bar x_p(y)}.
\end{equation}
If $\eta$ is regular enough, then the assumptions on $\beta_p$ and
$\bar x_p$ \eqref{eqAssumptions}--\eqref{assumptions:nonlinear_xp}
are fulfilled for all $r$, and the numerical scheme we propose
applies. In order to prove the convergence towards the solution of
\eqref{eq1} with the parameter function $\eta$, we observe that
the distance between $\eta$ and its approximation \eqref{mmm}
expressed in terms of the proper norm can be bounded by $C/r$,
where $C$ does not depend on $t, \mu$ and $y$ due to
\eqref{mmm2}--\eqref{mmm}. Thus, the most general version of the
stability result in {\rm\cite[Theorem 2.11]{Spop}} guarantees that
if $r$ tends to $+ \infty$, then the numerical solution obtained
for the approximated $\eta$ converges towards a solution to
\eqref{eq1} with the parameter function $\eta$. For all technical
details, we refer to {\rm\cite{Spop}}.
\end{remark}


\section{Simulation Results}\label{simulations}

This section is devoted to presenting results of numerical
simulations for several test cases. In all examples presented in
this paper, we used the $4$-th order Runge-Kutta method for
solving \eqref{b} and the explicit Euler scheme for solving
\eqref{c_eta}, as described in Subsection \ref{GD}. The error of
the numerical solution with parameters $(\Delta t, \bar M_o, \bar
M)$ at time $T>0$ is defined as
\begin{equation}\label{error:def}
\mathrm{Err}(T; \Delta t, \bar M_o, \bar M) :=
\rho(\mu(t_{\bar{k}}), \mu_{\bar{k}} )\,,
\end{equation}
with $\bar{k}$ such that $\bar{k}\Delta t=T$. The order of the
method $q$ is given by
\begin{equation}
q := \lim_{\Delta t \to 0}\frac{\log \left (\mathrm{Err}{(T;
2\Delta t, 2\bar M_o, 2\bar M)} / \mathrm{Err}{(T; \Delta t, \bar
M_o, \bar M)} \right)}{\log 2}.
\end{equation}
We also define $\Delta x := \modulo{K}/ \bar M_o$, where $K$ is
the minimal bounded closed interval containing the support of the
initial measure. We will not distinguish between measures and
their densities whenever the measures are absolutely continuous
with respect to the Lebesgue measure.

\subsection{Example 1 (McKendrick-von Foerster equation)}
In this subsection, we validate the convergence result for our
numerical scheme by means of the well-known McKendrick-von
Foerster type equation \cite{McK}. This is a linear model
describing the evolution of an size-structured population. We set
\begin{eqnarray*}
b(x) = 0.2(1-x), \;\; c(x) = 0.2, \;\; \eta(y)= 2.4(y^2 -
y^3)\delta_{x = 0} , \;\; \mathrm{and} \;\; \mu_o =
\chi_{[0,1]}(x),
\end{eqnarray*}
and solve \eqref{eq1} for $x \in [0,1]$, see also \cite{angulo1}.
The solution is stationary and then given by $\mu(t,x) =
\chi_{[0,1]}(x)$. In Table \ref{error:ex1:T}, we present the
relative error and the order of the scheme, where we used just one
measure reconstruction in order to approximate the initial data.
In Table \ref{error:ex1:T2}, we present results for the scheme
with the measure reconstruction performed at $t = 0, 1, \dots, 10$
and $\bar M_o = \bar M$. In all cases, we see that the convergence
error approximates order one as $\Delta t\to 0$ as proven in
Theorem \ref{rec_error} and Remark \ref{rmerror}.

\begin{table}[ht]{\tiny
\begin{center}
\begin{tabular}{ r  > {$\quad\quad\quad} c < {$}  > {$} c  < {$} }
\hline
 $\Delta t$ = $\Delta x$ &
\mathrm{Err}(10, \Delta t, \bar M_o, \bar M) &  q
\\
\hline $1.000000 \cdot 10^{-1} $ &1.2532 \cdot10^{-2} & -
\\
\hline $5.000000 \cdot 10^{-2} $ &5.0543  \cdot10^{-3}&  1.31006
\\
\hline $2.500000 \cdot 10^{-2} $ &2.2225\cdot10^{-3} &  1.18533
\\
\hline $1.250000 \cdot 10^{-2}$ &1.0349\cdot10^{-4} &  1.10272
\\
\hline $6.250000 \cdot 10^{-3} $ & 4.9832\cdot10^{-4} &  1.05431
\\
\hline $3.125000 \cdot 10^{-3} $ & 2.4438\cdot10^{-4} &  1.02796
\\
\hline $1.562500 \cdot 10^{-3} $ & 1.2099\cdot10^{-4} &  1.01419
\\
\hline $7.812500 \cdot 10^{-5} $ & 6.0198\cdot10^{-5} &  1.00715
\\
\hline $3.906250 \cdot  10^{-4} $ & 3.0024\cdot10^{-5} & 1.00359
\\
\hline $1.953125 \cdot 10^{-4} $ & 1.4993\cdot10^{-5} & 1.00180
\\
\hline $9.765625 \cdot  10^{-5} $ & 7.4920\cdot10^{-6} & 1.00090
\\
\hline
\end{tabular}
\caption{(Example 1) The relative error and order of the scheme at
$T = 10$. One reconstruction performed at $t = 0$, $\bar M = \bar
M_o$.}\label{error:ex1:T}
\end{center}}
\end{table}

\begin{table}[ht]{\tiny
\begin{center}
\begin{tabular}{ r > {$\quad\quad} c < {$} > {$} c  < {$} > {$} c < {$} > {$} c  < {$} }
\hline
 $\Delta t$ = $\Delta x$ &
\mathrm{Err}(10, \Delta t, \bar M_o, \bar M)& q&\mathrm{Err}(10, \Delta t, \bar M_o, \bar M)& q
\\
  & \mbox{(Fixed-location)} &  & \mbox{(Fixed-equal mass)} &
\\
\hline
$1.000000 \cdot 10^{-1} $ & 3.4657 \cdot10^{-1} & - & 8.8838  \cdot10^{-2} & -
\\
\hline
$5.000000\cdot 10^{-2} $ &1.1670 \cdot10^{-1}& 1.5703 & 2.9437 \cdot  10^{-2} & 1.5935
\\
\hline
$2.500000\cdot 10^{-2} $ &3.4080\cdot10^{-2} & 1.7759 & 1.0879 \cdot  10^{-2} & 1.4361
\\
\hline
$1.250000\cdot 10^{-2} $ &1.1863\cdot10^{-2} & 1.5224 & 4.4725 \cdot  10^{-3} & 1.2824
\\
\hline
$6.250000\cdot 10^{-3} $ & 3.6874\cdot10^{-3} & 1.6858 & 1.9907 \cdot  10^{-3} & 1.1678
\\
\hline
$3.125000\cdot 10^{-3} $ & 1.6866\cdot10^{-3} & 1.1285 & 9.3351 \cdot  10^{-4} & 1.0926
\\
\hline
$1.562500\cdot 10^{-3} $ & 6.8067\cdot10^{-4} & 1.3091 & 4.5131 \cdot  10^{-4} & 1.0486
\\
\hline
$7.812500 \cdot 10^{-4} $ & 3.3212\cdot10^{-4} & 1.0352 & 2.2178 \cdot  10^{-4} & 1.0250
\\
\hline
$3.906250 \cdot  10^{-4} $ & 1.5814\cdot10^{-4} & 1.0705 & 1.0992 \cdot  10^{-4} & 1.0127
\\
\hline
$1.953125  \cdot 10^{-4} $ & 7.4507\cdot10^{-5} & 1.0858 & 5.4719 \cdot  10^{-5} & 1.0063
\\
\hline
$9.765625 \cdot 10^{-5} $ & 3.6414\cdot10^{-5} & 1.0329 & 2.7299 \cdot  10^{-5} & 1.0032
\\
\hline
\end{tabular}
\caption{(Example 1) The relative error and order of the scheme at
$T = 10$. Reconstruction performed at $t = 0, 1, \dots, T$, $\bar
M = \bar M_o$.} \label{error:ex1:T2}
\end{center}}
\end{table}

\subsection{Example 2 (nonlinear growth term)}
In this subsection, we present results for a model where $b$ and
$\eta$ are equal to zero. Thus, we have conservation of the number
of approximated Dirac Deltas, and consequently, there is no need
for reconstructions. We consider a nonlinear growth function $c$
as in \cite{des} of the form
\begin{eqnarray*}
c(t,\mu)(x) = a(x) - \int_{\reali}\alpha(x,y)\d \mu(y),
\end{eqnarray*}
where
\begin{displaymath}
a(x) = A - x^2,\;\;A > 0
\;\;\;\;\; \mathrm{and} \;\;\;\;\;
\alpha(x,y) = \frac{1}{1 + (x-y)^2}.
\end{displaymath}
According to \cite[Remark 2.3, Lemma 4.8]{Spop}, one can consider
\eqref{eq1} on the whole $\reali$, so that the result concerning
well posedness still holds if all parameter functions verify the
regularity properties
\eqref{eqAssumptions}--\eqref{assumptions:nonlinear_xp} on the
whole line. However, $a(x)$ is not globally Lipschitz on $\reali$.
Nevertheless, the global well-posedness theory still applies if we
reduce to measures whose support lies in a fixed compact interval.
Note that the support of the solution is invariant in time.

\begin{figure}[h]
\begin{center}
\includegraphics[width=330px]{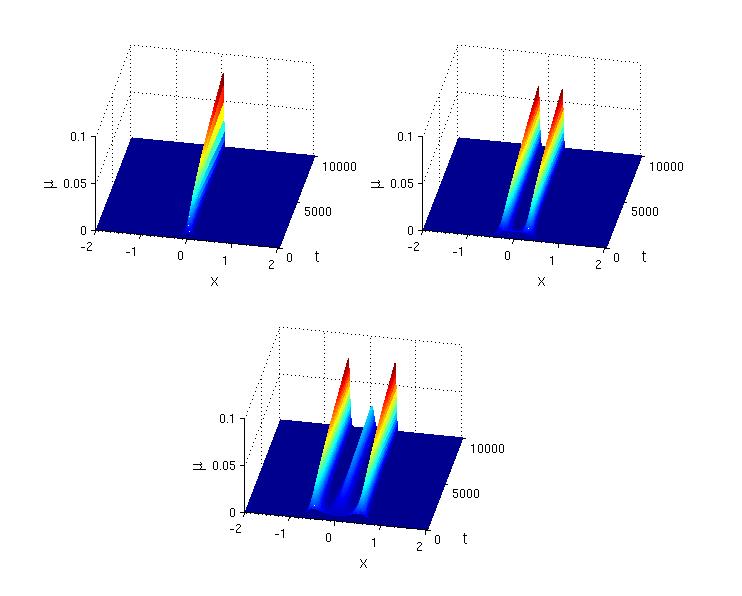}
\caption{(Example 2) Long time behaviour of numerical solutions.
The three subplots show the evolution of the numerical solution on
the time interval $[0, 10000]$ for $A = 0.5, 1.5$ and $2.5$,
respectively. For simulations, we set $\Delta t = 0.1$, $\bar M_o
= 1000$ and $\mu_o = \sum_{i=1}^{\bar M_o} ({1}/{\bar M_o})
\delta_{x^i_o}$, where $x^i_o := -2 + (i-\frac{1}{2})/\bar M_o$.
No measure reconstruction has been performed.} \label{figure2}
\end{center}
\end{figure}

\begin{figure}[h]
\begin{center}
\includegraphics[width=220px]{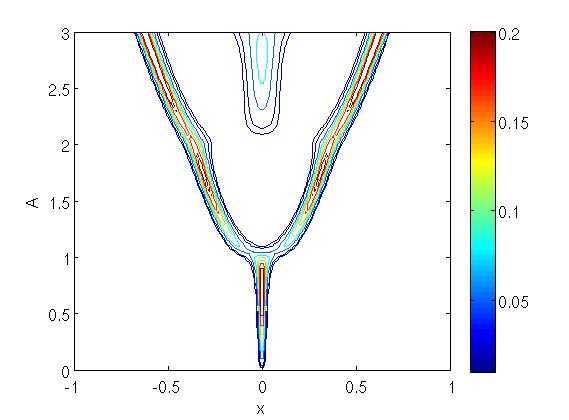}
\caption{(Example 2) Stationary State as a function of $A>0$. We
show the numerical solution at time $t=10000$ depending on the
parameter $A \in [0,3]$. For simulations, we set $\Delta t =
0.05$, $\bar M_o = 320$ and $\mu_o = \sum_{i=1}^{\bar M_o}
({1}/{\bar M_o}) \delta_{x^i_o}$, where $x^i_o := -2 +
(i-\frac{1}{2})/\bar M_o$. No measure reconstruction has been
performed.} \label{figure2b}
\end{center}
\end{figure}

If $\modulo{x} > \sqrt{A}$, then the solution decreases
exponentially to zero, since $\alpha(x,y) \geq 0$, for all $x,y
\in \reali$. This equation can describe a population structured
with respect to the trait $x$, and then its asymptotic behaviour
reflects the speciation process. Typically, after a long time
period only a few traits are observable, since the rest of the
population got extinct. Under some assumptions, there exists a
linearly stable steady solution $\bar \mu$ being a sum of Dirac
Deltas, which is shown in \cite{des}. The number of Dirac measures
depends on the parameter $A$ and some stationary solutions are
explicit. Figures \ref{figure2} and \ref{figure2b} present the
evolution and long time behaviour of solutions for different
choices of the parameter $A$. These results are consistent with
the findings in \cite{des}. In all cases, we assumed that initial
data are given as a sum of uniformly distributed Dirac Deltas with
the same mass.

\subsection{Example 3 (size structure - equal fission)}

In this subsection, we shall concentrate on a size-structured cell
population model, in which a cell reproduces itself by fission
into two equal parts. We assume that the cell divides after it has
reached a minimal size $x_o > 0$. Therefore, there exists a
minimum size whose value is $x_o / 2$. Moreover, cells have to
divide before they reach a maximal size, which is normalized to be
equal to $x_{max} = 1$. Similarly as in \cite{AnguloFission}, we
set
$$
x_o = \frac{1}{4}, \;\; b(x) = 0.1(1-x), \;\; c(x) =
\beta(x),\;\;\eta(t,\mu)(y) = 2\beta(y)\delta_{x =
y/2},\;\;\mathrm{and}\;\; \mu_o(x) = (1- x)(x - x_o/2)^3,
$$
where
\[
   \beta(y) = \left\{
   \begin{array}{c l}
    0, & \quad  \text{for $y \in (\reali^{+} \backslash\;
    [x_o,1])$,}\\[3mm]
    \displaystyle\frac{b(y) \phi(y)}{1 - \int_{x_o}^{y} \phi(x) \d x}, & \quad \text{for $y \in [x_o,1]$,}\\
   \end{array} \right.
 \]
and
\[
   \phi(y) = \left\{
   \begin{array}{l l}
\frac{160}{117}\left( -\frac{2}{3} + \frac{8}{3} y \right)^3,
& \quad  \text{for $y \in  [x_o,(x_o + 1)/2]$,}\\
\frac{32}{117}\left( -20 + 40 y + \frac{320}{3}\left(y -
\frac{5}{8} \right)^2 \right) + \frac{5120}{9}\left(y -
\frac{5}{8} \right)^3\left( \frac{8}{3} y - \frac{11}{3} \right),
 & \quad \text{for $y \in ((x_o +1)/2,1]$.}\\
   \end{array} \right.
 \]

\begin{figure}[h]
\begin{center}
\includegraphics[width=300px]{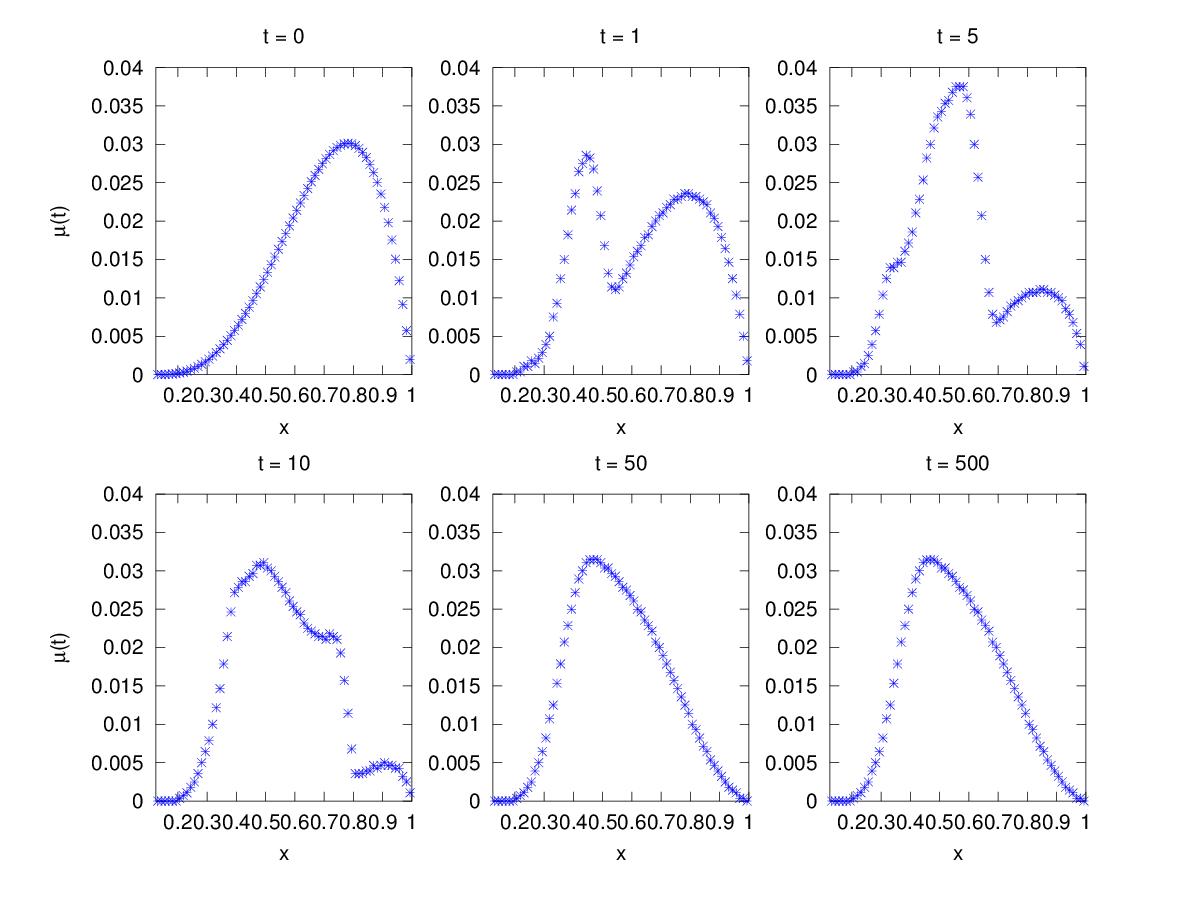}
\caption{(Example 3) Numerical solution at $t = 0, 1, 5, 10, 50,
500$, calculated for $\Delta t = 0.0125$, $\bar M_o = \bar M =
2800$. Fixed-equal mass reconstruction has been performed every
$4$ time steps. We show the numerical solution after the
fixed-location reconstruction with parameter $\bar M = 70$ and
normalization.} \label{figure23}
\end{center}
\end{figure}

Figure \ref{figure23} shows the long time behaviour of a numerical
solution for a particular choice of parameters. We observe the
convergence towards a stationary profile once normalized, since
the mass grows exponentially in time, as discussed in
\cite{DHT,AnguloFission}. We remark that this structured
population model cannot be discretized using the standard EBT
method since particles divide at different sizes and the nonlocal
term cannot be understood as a boundary condition. In order to
keep the number of Dirac Deltas under control, we perform the
reconstruction procedure as discussed in Subsection 2.3. Let us
point out that the convergence towards normalized stationary
states for similar models in the framework of Lebesgue spaces has
been proved in \cite{PR,LP,CCM,BCG}. Finding the properties of
these stationary states numerically is a relevant question that
will be discussed elsewhere.

\subsection{Example 4  (selection-mutation)}

The last test case concerns a simple selection-mutation model in
which the population is structured with respect to a evolutionary
trait as in \cite{CCDR}. We assume that $x \in [0,1]$ and set the
parameters as
$$
b(x) = 0, \;\; c(\mu)(x) = - (1-\epsilon)B(x) + m(\mu),\;\;\;
\mathrm{and}\;\;\; \eta(y) = \epsilon \sum_{p=1}^{r}B(y)
\beta_p(y)\delta_{x = \bar x_p(y)}.
$$
Here, $B(x)$ represents the trait specific birth rate, $m(\mu)$ is
the death rate depending on the population distribution, and
$\beta_p$ represents the mutation density probability, i.e., the
probability that a parent with trait $y$ has a newborn with trait
$\bar x_p(y)$. Finally, the parameter $\epsilon$ is the mutation
rate, and thus there are two parts in the right hand side, those
that are a faithful reproduction of their parents and those that
mutate, slightly with high probability, their trait.

Let us point out that the mutation term in this model is an
approximation in the sense of Remark \ref{aproxeta} of a
continuous nonlocal term of the form
$$
\int_0^1 B(y) \beta(x,y) \,\d\mu(y)\,\qquad \mbox{with}\qquad
\int_0^1 \beta(x,y) \d x = 1\,,
$$
and, in practice we can assume that has a Gaussian shape
concentrated around the diagonal $x=y$. The approximated nonlocal
term is constructed by substituting the mutation probability
density $\beta(x,y)$ at each $y$ by an approximation with $r$
Delta Dirac points $\{\bar x_p(y)\}_{p=1}^r$ leading to the form
of $\eta(y)$ above. More precisely, the approximated $\eta(y)$ is
defined by duality on test functions $\varphi\in \mathbf{C}_0(\reali^+)$
functions as
\begin{align*}
\int_{\reali^+} \int_{\reali^+}\varphi(t,x) B(y) \beta(x,y) \d x
\d \mu(y) &\approx \sum_{p=1}^{r} \int_{\reali^+} B(y) \beta_p(y)
\varphi(\bar x_p(y)) \d\mu(y) \\
&= \int_{\reali^+} \int_{\reali^+}\varphi(t,x)[\d\eta(t,\mu)(y)](x)
\d\mu(y) \,.
\end{align*}

\begin{figure}[h]
\begin{center}
\includegraphics[width=230px]{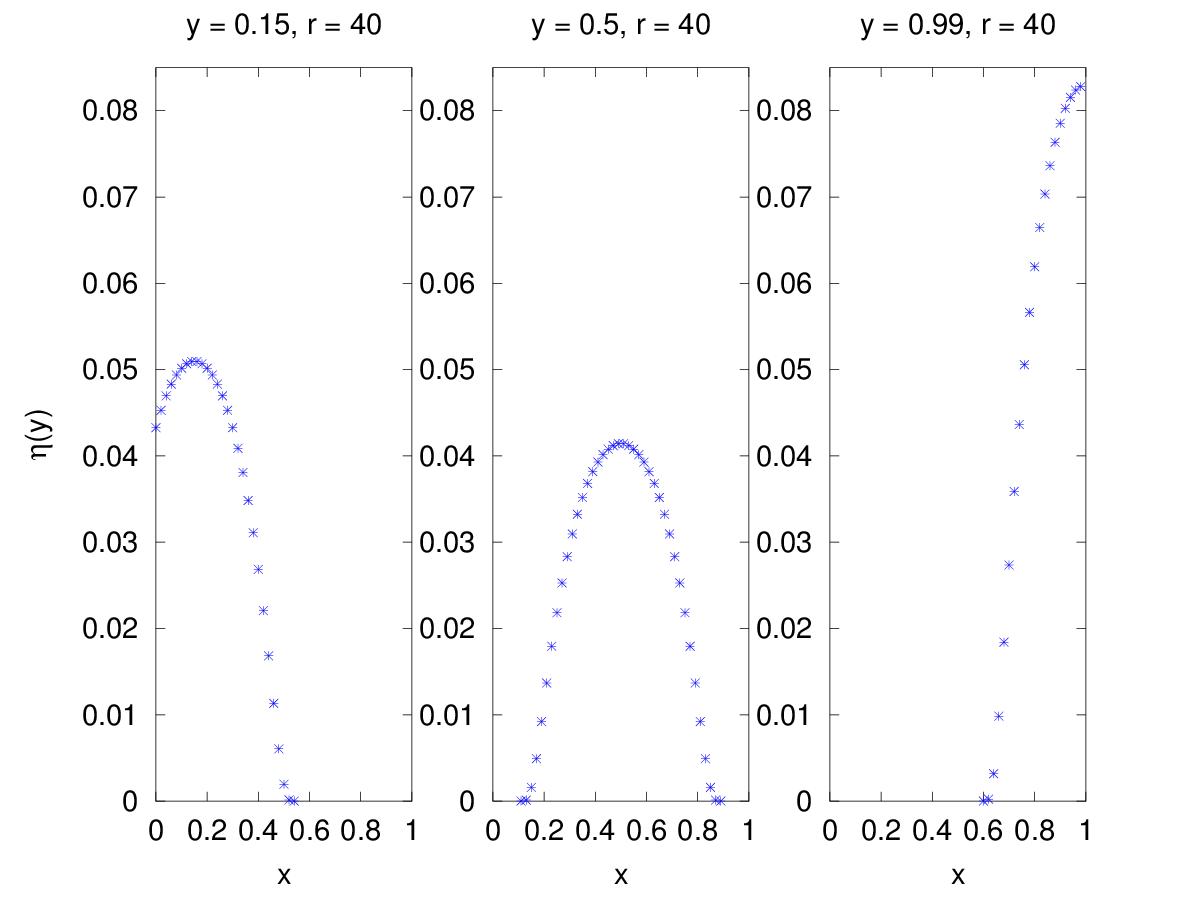}
\caption{(Example 4) The subplots show the function $\eta(y)$ for
$y = 0.15$, $y=0.5$ and $y=0.99$, respectively, and parameters $r
= 40$, $a = 0.4$. } \label{eta}
\end{center}
\end{figure}

In our simulations and based on the previous considerations, we
consider $B(x) = x(1-x)$, the death rate is assumed to depend
increasingly on the total population with a saturation of the form
$m(\mu) = 1 - \exp\left\{-\int_{0}^{1} \d \mu\right\}$, and the
approximation of the mutation kernel is chosen with $r = 10$,
\[
   \bar x_p(y) = \left\{
   \begin{array}{l l}
     (y-a) + \frac{a}{r}\left(2p - 1\right), & \quad  \text{if $0 \leq (y-a) + \frac{a}{r}\left(2p - 1\right) \leq 1$,}\\
    0, & \quad \text{otherwise},\\
   \end{array} \right.
 \]
and
\[
\beta_p(y) = \frac{\check \beta_p(y)}{\sum_{p=1}^{r} {\check \beta_p(y)}},\;\;
\text{where}\;\;
   \check \beta_p(y) = \left\{
   \begin{array}{l l}
    \mathrm{exp}\left(
-\frac{a^2}{a^2 - (\bar x_p(y) - y)^2}
\right)
, & \quad  \text{if $p$ is s.t. $ 0 \leq \bar x_p(y) \leq 1$,}\\
    0, & \quad \text{otherwise}.\\
   \end{array}  \right..
 \]
The parameter $a$ is related to the mutation strength in the sense
that a distance between a parent and its offspring is not greater
than $a$, set in our simulations to $a = 0.4$.

\begin{figure}[h]
\begin{center}
\includegraphics[width=300px]{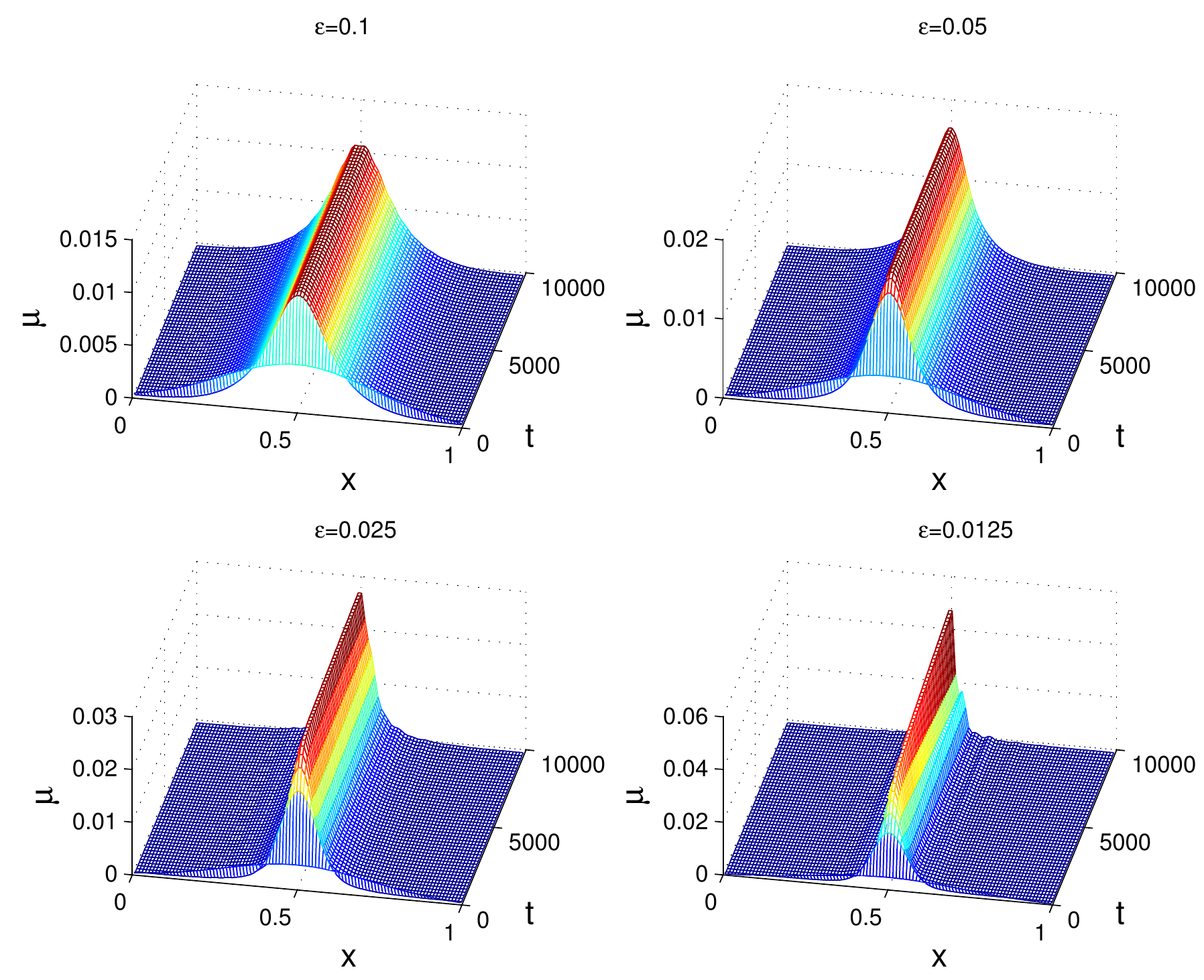}
\caption{(Example 4) Long time behaviour of numerical solutions.
The plots show the evolution of a numerical solution on the time
interval $[0, 2000]$ for $\epsilon = 0.1, 0.05, 0.025$, and
$0.0125$, respectively. For simulations, we set $\Delta t =
0.025$, $\bar M_o = \bar M= 100$, and $\mu_o = \sum_{i=1}^{\bar
M_o} ({1}/{\bar M_o}) \delta_{x^i_o}$, where $x^i_o :=
(i-\frac{1}{2})/\bar M_o$. Fixed location reconstruction has been
performed every $2$ time steps.} \label{ex4_1}
\end{center}
\end{figure}

Figure \ref{ex4_1} shows the convergence towards stationary states
for different values of the mutation rate $\epsilon$. We observe
that the stabilization rate depends on $\epsilon$, being slower as
$\epsilon$ gets smaller and smaller. The existence of these
stationary states with the full mutation kernel $\eta$ was proved
in \cite{CCDR} without information about their stability.

\section*{Acknowledgments}
JAC acknowledges support from the Royal Society by a Wolfson
Research Merit Award and by the Engineering and Physical Sciences
Research Council grant with references EP/K008404/1. JAC was
partially supported by the project MTM2011-27739-C04-02 DGI
(Spain) and 2009-SGR-345 from AGAUR-Generalitat de Catalunya.
PG
is the coordinator and AU is a Ph.D student in the International
Ph.D. Projects Programme of Foundation for Polish Science operated
within the Innovative Economy Operational Programme 2007-2013
(Ph.D. Programme: Mathematical Methods in Natural Sciences). PG is
supported by the grant of National Science Centre no
6085/B/H03/2011/40. AU is supported by the grant of National Science Centre no 2012/05/N/ST1/03132.

\bibliographystyle{plain}
\bibliography{Particle}

\end{document}